\documentclass{amsart}

\newtheorem{theorem}{Theorem}[section]
\newtheorem{lemma}[theorem]{Lemma}

\theoremstyle{definition}

\newtheorem{cor}[theorem]{Corollary}
\usepackage{amsmath}
\usepackage{centernot}
\usepackage{amssymb,latexsym}
\usepackage{enumerate}
\usepackage{cases}
\usepackage{mathtools}
\theoremstyle{remark}

\numberwithin{equation}{section}

\begin{document}

\baselineskip=16pt

\title[Double square moments and bounds for resonance sums]{Double square moments and bounds for
\\
resonance sums of cusp forms}

\author{Tim Gillespie, Praneel Samanta, and Yangbo Ye}

\address{Tim Gillespie${}^1$: GillespieTimothyL@sau.edu}

\address{Praneel Samanta${}^2$: praneel-samanta@uiowa.edu}

\address{Yangbo Ye${}^2$: yangbo-ye@uiowa.edu}

\address{${}^1$
Department of Mathematics and Statistics, St. Ambrose 
University, Davenport, Iowa 52803, USA }

\address{${}^2$ Department of Mathematics, The University of 
Iowa, Iowa City, Iowa 52242, USA}

\thanks{These authors have contributed equally to this work and share first authorship.}



\begin{abstract}
Let $f$ and $g$ be holomorphic cusp forms for the modular 
group $SL_2(\mathbb Z)$ of weight $k_1$ and $k_2$ with 
Fourier coefficients $\lambda_f(n)$ and $\lambda_g(n)$, 
respectively. For real $\alpha\neq0$ and 
$0<\beta\leq1$, consider a smooth resonance sum 
$S_X(f,g;\alpha,\beta)$ of $\lambda_f(n)\lambda_g(n)$ 
against $e(\alpha n^\beta)$ over $X\leq n\leq2X$. 
Double square moments of $S_X(f,g;\alpha,\beta)$ 
over both $f$ and $g$ are nontrivially bounded 
when their weights $k_1$ and $k_2$ tend to infinity 
together. By allowing both $f$ and $g$ to move, these 
double moments are indeed square moments associated with 
automorphic forms for $GL(4)$. By taking out a small 
exceptional set of $f$ and $g$, bounds for individual 
$S_X(f,g;\alpha,\beta)$ will then be proved. These 
individual bounds 
break the resonance barrier of $X^\frac58$ for 
$\frac16<\beta<1$ and achieve a square-root cancellation 
for $\frac13<\beta<1$ for almost all $f$ and $g$ as an 
evidence for Hypothesis S for cusp forms over integers. 
The methods used in this study include Petersson's 
formula, Poisson's summation formula, and stationary 
phase integrals. 
\end{abstract}

\maketitle

\section{Introduction} 

According to Iwaniec, Luo, and Sarnak \cite[Appendix C]{ILS}, 
a general form of Hypothesis S over integers predicts a 
square-root cancellation in the sum
$$
S_X(\{a_n\};\alpha,\beta)
=
\sum_na_n
e(\alpha n^\beta)
\phi\Big(\frac nX\Big)
\ll
X^{\frac12+\varepsilon},
$$
where $0<\beta\leq1$ and $\alpha \neq 0$ are fixed real numbers, 
$\phi$ is a smooth function of compact support in $(1,2)$, and 
$\{a_n\}$ is an arithmetically defined sequence of complex numbers 
satisfying $a_n\ll n^\varepsilon$. The Case III in \cite{ILS} of 
$a_n=\lambda_f(n)$ for a fixed holomorphic cusp form $f$ for 
$SL_2(\mathbb Z)$, however, faces a resonance barrier when 
$\beta=\frac12$ and $\alpha=\pm2\sqrt q$ for a positive integer 
$q$. In other words, 
$$
S_X\Big(f;\pm2\sqrt q,\frac12\Big)
:=
\sum_n 
\lambda_f(n)
e(\pm2\sqrt{nq})
\phi\Big(\frac nX\Big)
$$ 
has a main term of size $|\lambda_f(q)|X^\frac34$. This resonance 
phenomenon has been further studied by Ren and Ye 
\cite{RY10} -- \cite{RY16}, Ernvall-Hyt\"onen \cite{EH}, 
Ernvall-Hyt\"onen, J\"a\"asaari, and Vesalainen \cite{EHJV}, 
Czarnecki \cite{KC}, Savala \cite{Sava}, and many other authors 
for fixed automorphic forms $f$.

It is believed that one might be able to break the resonance 
barrier if the cusp form $f$ is allowed to move. In this direction, 
Ye \cite{Ye22} proved the first known non-trivial bound for 
$S_X(f;\alpha,\beta)$ when the weight $k$ of $f$ tends to infinity with the summation length $X$. 
This bound, however, is far from reaching the square-root 
cancellation.

The goal of the present paper is to break the resonance barrier 
for  
\begin{equation}\label{S_X}
S_X(f, g;\alpha,\beta)
=
\sum_n
\lambda_f(n)
\lambda_g(n)
e(\alpha n^{\beta})
\phi\Big(\frac{n}{X}\Big)
\end{equation}
for almost all holomorphic cusp forms $f$ and $g$ for 
$SL_2(\mathbb{Z})$ of even integer weights $k_1$ and $k_2$, 
respectively. Note here that 
$\lambda_f(n)\lambda_g(n)$ in \eqref{S_X} 
corresponds to the Dirichlet coefficients coming from the 
Rankin-Selberg $L$-function 
$$
L(s,f\times g) 
= 
\zeta(2s)\sum_{n\geq 1} \frac{\lambda_f(n)\lambda_g(n)}{n^s}.      
$$
Consequently, $S_X(f,g;\alpha,\beta)$ represents the interplay 
between the oscillation of the Dirichlet coefficients of 
$L(s,f\times g)$ and that of a fractional exponential 
function. 

By \cite{RY15SCM} and \cite{KC}, the sum 
$S_X(f,g;\alpha, \beta)$ has a main term of size 
$|\lambda_f(q)\lambda_g(q)|X^{\frac{1}{8}+\frac{1}{2}}$ 
when $\beta=\frac14$ and $\alpha$ is close or equal to 
$\pm 4q^\frac14$ for a positive integer $q$ for fixed 
$f$ and $g$. The resonance 
barrier in this case is thus $X^\frac58$. We will break 
this resonance barrier for almost all $f$ and $g$, in the 
sense that the number of exceptional pairs $(f,g)$ 
is a power less than the total number of pairs $(f,g)$ in 
consideration. In particular, the exceptional pairs have a 
probability tending to zero. 

More precisely, let $S_k(\Gamma)$ denotes the space of 
holomorphic cusp 
forms on $\Gamma=SL_2(\mathbb{Z})$ of even integral weight 
$k$, and let $H_k$ denote an orthogonal basis of 
$S_k(\Gamma)$ consisting of Hecke eigenforms where each 
form is normalized to have the first Fourier coefficient 
equal to $1$.  Recall the dimension formula 
$\dim S_k(\Gamma) = \frac k{12}+O(1)$. Consider parameters 
$K_j^{\varepsilon }\leq L_j\leq K_j^{1-\varepsilon}$ for 
$j=1,2$, and denote 
$H_{K_j,L_j}=\cup_{K_j-L_j\leq k\leq K_j+L_j}H_k$. We will 
bound the double square moment 
\begin{equation}\label{SSumK1K2}
\sum_{f\in H_{K_1,L_1}} 
\sum_{g\in H_{K_2,L_2}} 
|S_X(f,g;\alpha,\beta)|^2
\end{equation}
for $K_1,L_1,K_2,L_2,X$ in various ranges. Let 
$g_0\in C^\infty(-1,1)$ be a non-negative test function 
with $g_0(0)=1$. We will equivalently bound the smooth 
double moment 
\begin{equation}\label{SumK1K2}
\sum_{2 | k_1}
\sum_{2 | k_2}
g_0\Big(\frac{k_1-K_1}{L_1}\Big)
g_0\Big(\frac{k_2-K_2}{L_2}\Big)
\sum_{f\in H_{k_1}}
\sum_{g\in H_{k_2}}
|S_X(f,g,\alpha,\beta)|^2.
\end{equation}
In order to apply Petersson's formula, we will actually 
bound a normalized sum
\begin{eqnarray}
\sum\nolimits_{K_1L_1}^{K_2L_2} 
&=&  
K_1K_2 
\sum_{2 | k_1}
\sum_{2 | k_2} 
g_0\Big(  \frac{k_1-K_1}{L_1}   \Big) 
g_0\Big(  \frac{k_2-K_2}{L_2}   \Big) 
\sum_{f\in H_{k_1} } 
\frac{2 \pi^2}{ (k_1-1)L(1, Sym^2 f)}
\label{NSumK1K2}
\\ 
&&\times 
\sum_{g\in H_{k_2} } 
\frac{2 \pi^2}{ (k_2-1)L(1, Sym^2 g)} 
|S_X(f, g,\alpha, \beta)|^2 
\nonumber 
\\ 
&=&    
K_1K_2 
\sum_{2 | k_1}
\sum_{2 | k_2} 
g_0\Big(  \frac{k_1-K_1}{L_1}   \Big) 
g_0\Big(  \frac{k_2-K_2}{L_2}   \Big) 
\sum_n 
\sum_m 
e(\alpha n^{\beta}-\alpha m^{\beta}) 
\phi\Big(\frac{n}{X}\Big)
\bar{\phi}\Big(\frac{m}{X}\Big)
\nonumber
\\
&&\times
\sum_{f\in H_{k_1} }
\frac{2 \pi^2\lambda_f(n)\bar{\lambda}_f(m)}{ (k_1-1)L(1, Sym^2 f)}
\sum_{g\in H_{k_2} } 
\frac{2 \pi^2\lambda_g(n)\bar{\lambda}_g(m)}{ (k_2-1)L(1, Sym^2 g)}.
\nonumber
\end{eqnarray}
\bigskip
The factors 
$$
\frac{2\pi^2K_1}{(k_1-1)L(1,Sym^2 f)}
\ \text{and}\ 
\frac{2\pi^2K_2}{(k_2-1)L(1,Sym^2 g)}
$$ 
in \eqref{NSumK1K2} will result 
in a discrepancy bounded between $K_i^{-\varepsilon}$ and 
$K_i^{\varepsilon}$, $i=1,2$, by Iwaniec \cite{I} and 
Hoffstein and Lockhardt \cite{HL}. 

For holomorphic cusp forms we have Deligne's proof \cite{D} of 
the Ramanujan conjecture which implies that 
$\lambda_f(n)\ll n^\varepsilon$ and 
$\lambda_g(n)\ll n^\varepsilon$ for $\varepsilon>0$ 
arbitrarily small. Here the implied constants are independent of 
$f$ and $g$. Using this we obtain a trivial bound 
\begin{equation}
\label{TrivialB}
S_X(f,g;\alpha,\beta)
\ll
X^{1+\varepsilon},
\end{equation} 
where the implied constant is independent of $f$ and $g$. 
Applying this to \eqref{NSumK1K2}, we get a trivial bound 
$O(K_1L_1K_2L_2 X^{2+\varepsilon})$ 
for $\sum\nolimits_{K_1L_1}^{K_2L_2}$ in \eqref{NSumK1K2} 
and hence for \eqref{SumK1K2} and \eqref{SSumK1K2}. We 
seek to break this bound. Note that 
non-trivial bounds beyond \eqref{TrivialB} are known for 
$S_X(f,g;\alpha,\beta)$ for both $f$ and $g$ being fixed 
but not for $f$ and $g$ with their 
weights tending to infinity.

\begin{theorem}\label{Thm1}
For $j=1,2$ assume $K_j^{\varepsilon } \leq L_j \leq K_j^{1-\varepsilon}$. Then for $0<\beta<1$
\begin{eqnarray}
\sum\nolimits_{K_1L_1}^{K_2L_2} 
&\ll&
K_1L_1K_2L_2X^{1+\varepsilon}
\ \ {\rm if}\,
K_1L_1\geq X^{1+\varepsilon}
\,{\rm and}\,
K_2\geq X^{\frac12+\varepsilon};
\label{K1K2>}
\\
&\ll&
K_1L_1K_2L_2X^{1+\varepsilon}
+
\frac{K_1L_1L_2}{K_2}X^{\frac32+\varepsilon}
\ \ {\rm if}\,
K_1L_1\geq X^{1+\varepsilon}
\,{\rm and}\,
K_2\leq X^\frac12;
\label{K1>K2<}
\\
&\ll&
K_1^2L_1L_2X^{1+\varepsilon}
+
\frac{X^{3+\varepsilon}}{K_1}
\label{K1=K2>}
\\
&&
{\rm if}\,
K_1L_1,K_2L_2\leq X^{1+\varepsilon},\,
K_1=K_2,
\,{\rm and}\,
K_1^2L_1L_2\geq X^{1+\beta+\varepsilon}.
\nonumber
\end{eqnarray}
When $\beta=1$, bounds in \eqref{K1K2>} and 
\eqref{K1>K2<} remain valid, while \eqref{K1=K2>} 
is replaced by 
\begin{eqnarray}\label{K1=K2beta=1}
&\ll&
\min(L_1,L_2)K_1X^{2+\varepsilon}
+
\frac{X^{3+\varepsilon}}{K_1}
\ \ {\rm if}\,
K_1L_1,K_2L_2\leq X^{1-\varepsilon},\,
K_1=K_2.
\end{eqnarray}
\end{theorem}

Since the number of terms in 
$\sum\nolimits_{K_1L_1}^{K_2L_2}$ 
is $\asymp K_1L_1K_2L_2$, \eqref{K1K2>} show that the 
average size of $S_X(f,g,\alpha,\beta)$ is bounded by 
$O(X^{\frac12+\varepsilon})$ when $K_1L_1$ and $K_2$ 
are large. This represents a square-root saving on average 
for $S_X(f,g,\alpha,\beta)$. 

If we take 
\begin{equation}\label{KLrange}
K_1=K_2=X^\frac{1+\beta}2,
\ \ \  
L_1=L_2=X^\frac\varepsilon2
\end{equation}
in \eqref{K1=K2>}, we have 
$K_1^2L_1L_2=X^{1+\beta+\varepsilon}$ and the bound in 
\eqref{K1=K2>} becomes 
\begin{equation}\label{BonBeta}
\sum\nolimits_{K_1L_1}^{K_2L_2} 
\ll
X^{2+\beta+\varepsilon}
+
X^{\frac{5-\beta}2+\varepsilon}.
\end{equation}
This is also a bound for \eqref{SSumK1K2} for $0<\beta<1$. 
Since the two terms on the right hand side of \eqref{BonBeta} 
are both larger than $X^2$, \eqref{BonBeta} cannot provide a 
non-trivial bound for all individual $S_X(f,g;\alpha,\beta)$. 
It can, nevertheless, allow us to break the resonance barrier 
$X^\frac58$ for almost all forms $f$ and $g$. 

To this end, define 
$$
\Delta_{f,g}(\gamma)
=
\Big\{
{\rm pairs}\ (f,g)
\Big|
f\in H_{K_1,L_1}
,\,
g\in H_{K_2,L_2}
\,{\rm such\,that}\,
\big|S_X(f,g;\alpha,\beta)\big|
\geq X^{\gamma+\varepsilon}
\Big\}
$$
for $0<\gamma<1$. By the bound \eqref{BonBeta} on 
\eqref{SSumK1K2} we get 
$$
X^{2\gamma+2\varepsilon}
\big|\Delta_{f,g}(\gamma)\big|
\ll
X^{2+\beta+\varepsilon}
+
X^{\frac{5-\beta}2+\varepsilon}.
$$
Consequently, 
\begin{equation}\label{BDelta}
\big|\Delta_{f,g}(\gamma)\big|
\ll
X^{2+\beta-2\gamma-\varepsilon}
+
X^{\frac{5-\beta}2-2\gamma-\varepsilon}.
\end{equation}
To make the exceptional set $\Delta_{f,g}(\gamma)$ 
a power smaller than the size 
$K_1^2L_1L_2=X^{1+\beta+\varepsilon}$ of the averaging 
on $H_{K_1,L_1}$ and $H_{K_2,L_2}$, we require 
the two terms on the right hand side of \eqref{BDelta} 
$\leq X^{1+\beta}$. In other words, we need 
$2+\beta-2\gamma\leq1+\beta$ and 
$\frac{5-\beta}2-2\gamma\leq1+\beta$, i.e., 
\begin{equation}\label{gamma}
\max\Big(\frac12,\frac34(1-\beta)\Big)\leq\gamma<1. 
\end{equation}

\begin{cor}\label{IndividualB}
For parameters in \eqref{KLrange} and $\gamma$ in 
\eqref{gamma}, we have 
\begin{equation}\label{BreakB}
\big|S_X(f,g;\alpha,\beta)\big|
\leq
X^{\gamma+\varepsilon}
\end{equation}
for all $f\in H_{K_1,L_1}$ and $g\in H_{K_2,L_2}$ 
except for at most 
$O(X^{2+\beta-2\gamma}+X^{\frac{5-\beta}2-2\gamma})$ 
pairs of $f$ and $g$. The exceptional set is a power 
smaller than the size of $H_{K_1,L_1}\times H_{K_2,L_2}$.
\end{cor}

When we take 
$\gamma=\max\big(\frac12,\frac34(1-\beta)\big)$, the 
bound \eqref{BreakB} breaks the resonance barrier 
$X^\frac58$ when $\frac16<\beta<1$ and reaches the 
square-root cancellation $X^{\frac12+\varepsilon}$ 
when $\frac13<\beta<1$ for almost all $f$ and $g$. 
Similar corollaries can be formulated from 
\eqref{K1K2>}, \eqref{K1>K2<}, amd \eqref{K1=K2beta=1}. 

\setcounter{equation}{0}
\section{Petersson's Trace Formula} 

We recall Petersson's formula 
(Petersson \cite{P}, cf. Liu and Ye \cite{LiuYe}) for 
$m,n\geq 1$, 
$$     
\sum_{f\in H_k}    
\frac{2 \pi^2}{ (k-1)L(1, Sym^2 f)}
\lambda_f(n) 
\overline{\lambda}_f(m)
= 
\delta(m,n) 
+ 
2\pi i^k 
\sum_{c\geq 1} 
\frac{S(m,n,c)}{c}
J_{k-1}\Big( \frac{4\pi \sqrt{mn}}{c}   \Big),
$$
where $\lambda_f(n)$ is the Fourier coefficient of $f$ and 
$S(m,n,c)$ is the classical Kloosterman sum. Applying it to 
\eqref{NSumK1K2} we get
\begin{eqnarray}\label{4Terms}
\sum\nolimits_{K_1L_1}^{K_2L_2} 
& =&     
K_1K_2 
\sum_{2 | k_1}
\sum_{2 | k_2} 
g_0\Big(  \frac{k_1-K_1}{L_1}   \Big) 
g_0\Big(  \frac{k_2-K_2}{L_2}   \Big) 
\sum_n 
\sum_m 
e(\alpha n^{\beta}-\alpha m^{\beta}) 
\phi\Big(\frac{n}{X}\Big)
\bar{\phi}\Big(\frac{m}{X}\Big)
\\
&&\times
\Big(
\delta(n,m)
+
2\pi i^{k_1}
\sum_{c_1\geq 1}
\frac{S(m,n,c_1)}{c_1}
J_{k_1-1}\Big(\frac{4\pi\sqrt{mn}}{c_1}\Big) 
\Big) 
\nonumber
\\ 
&&\times
\Big(
\delta(n,m)
+
2\pi i^{k_2}
\sum_{c_2\geq 1}
\frac{S(m,n,c_2)}{c_2}
J_{k_2-1}\Big(\frac{4\pi\sqrt{mn}}{c_2}\Big) 
\Big)
\nonumber
\\
&=:& 
D_{00}+D_{01}+D_{10}+D_{11}.
\nonumber
\end{eqnarray}
\bigskip
Here
\begin{eqnarray}
D_{00} 
&=& 
K_1K_2
\sum_{2 | k_1}
\sum_{2 | k_2}
g_0\Big(\frac{k_1-K_1}{L_1}\Big)
g_0\Big(\frac{k_2-K_2}{L_2}\Big)
\sum_{n}
\Big|\phi\Big(\frac{n}{X}\Big)\Big|^2
\label{D00}
\\
D_{01} 
&=& 
K_1K_2 
\sum_{2 | k_1}
\sum_{2 | k_2} 
g_0\Big(  \frac{k_1-K_1}{L_1}   \Big) 
g_0\Big(  \frac{k_2-K_2}{L_2}   \Big) 
\sum_n 
\Big|\phi\Big(\frac{n}{X}\Big)\Big|^2
\label{D01}
\\
&&\times
2\pi i^{k_2}
\sum_{c_2\geq 1}
\frac{S(n,n,c_2)}{c_2}
J_{k_2-1}\Big(\frac{4\pi n}{c_2}\Big) 
\nonumber
\\
D_{10} 
&=& 
K_1K_2 
\sum_{2 | k_1}
\sum_{2 | k_2} 
g_0\Big(  \frac{k_1-K_1}{L_1}   \Big) 
g_0\Big(  \frac{k_2-K_2}{L_2}   \Big) 
\sum_n 
\Big|\phi\Big(\frac{n}{X}\Big)\Big|^2
\label{D10}
\\
&&\times
2\pi i^{k_1}
\sum_{c_1\geq 1}
\frac{S(n,n,c_1)}{c_1}
J_{k_1-1}\Big(\frac{4\pi n}{c_1}\Big) 
\nonumber
\end{eqnarray}
\begin{eqnarray}
D_{11} 
&=& 
K_1K_2 
\sum_{2 | k_1}
\sum_{2 | k_2} 
g_0\Big(  \frac{k_1-K_1}{L_1}   \Big) 
g_0\Big(  \frac{k_2-K_2}{L_2}   \Big) 
\sum_n 
\sum_m 
e(\alpha n^{\beta}-\alpha m^{\beta}) 
\phi\Big(\frac{n}{X}\Big)
\bar{\phi}\Big(\frac{n}{X}\Big)
\label{D11}
\\
&&\times
4\pi^2 i^{k_1+k_2}
\sum_{c_1\geq 1}
\frac{S(m,n,c_1)}{c_1}
J_{k_1-1}\Big(\frac{4\pi\sqrt{mn}}{c_1}\Big) 
\sum_{c_2\geq 1}
\frac{S(m,n,c_2)}{c_2}
J_{k_2-1}\Big(\frac{4\pi\sqrt{mn}}{c_2}\Big).
\nonumber
\end{eqnarray}

The diagonal term $D_{00}$ in \eqref{D00} can be 
estimated trivially for $0<\beta\leq1$: 
\begin{equation}\label{D00B}
D_{00}
\ll 
K_1L_1K_2L_2X.
\end{equation}

\setcounter{equation}{0}
\section{The $k_1$- and $k_2$-sums of Bessel functions} 

For off-diagonal terms $D_{01}$, $D_{10}$, and $D_{11}$, 
we have the quantity
$$
V_{K_j,L_j}(x) 
= 
\sum_{2 | k_j} 
i^{k_j} 
g_0\Big(\frac{k_j-K_j}{L_j}\Big)
J_{k_j-1}(x),
$$
with $x = \frac{4\pi \sqrt{mn}}{c_j}$.
Rewriting this sum as an oscillatory integral
$$  
V_{K_j,L_j}(x) 
= 
\frac{1}{2i}
\big( W_{K_j,L_j}(-x) - W_{K_j,L_j}(x),\big)
$$
where for $\eta=\pm 1$,
$$
W_{K_j,L_j}(\eta x) 
= 
\int_{-\infty}^{\infty} 
\hat{g}_0(t) 
e\Big( 
-\frac{(K_j-1)t}{L_j} 
- \frac{\eta x}{2\pi}\cos \frac{2\pi t}{L_j}         
\Big)dt,     
$$
then applying the method of stationary phase again gives 
an asymptotic expansion. Denote
$$
g(t)=\hat{g}_0(t),\ f(t)
=-\frac{(K_j-1)t}{L_j}-\frac{x}{2\pi} \cos \frac{2\pi t}{L_j}.
$$

\begin{lemma}\label{Wasymp}
{\rm (Sun and Ye \cite[Lemma 4.1]{SunYe}. cf. Salazar and Ye 
\cite{SalYe})} 
Suppose $K_j^\varepsilon\leq L_j\leq K_j^{1-\varepsilon}$. 

(1) If $x\leq 8\pi K_j^{1-\varepsilon}L_j$, then 
$W_{K_j,L_j}(\eta x)\ll_{\varepsilon, A}K_j^{-A}$ 
for any $A>0$.

(2) If $x\geq K_j^{1-\varepsilon}L_j$, then 
$$
W_{K_j,L_j}(\eta x) 
= 
\sum_{\nu = 0 }^{n_0}    
\widetilde{W}_{j,\nu}(\eta x) 
+ 
O\Big(  \frac{L_j^{2n_0+2}}{x^{n_0 + 1}}\Big).
$$
Here 
\begin{eqnarray*}   
\widetilde{W}_{j,\nu}(\eta x) 
&=&  
\frac{e(\frac18) (-1)^\nu (2\nu -1)!!}
{2^\nu i^\nu \eta^{\nu +\frac12} \pi^{2\nu +\frac12}}
\frac{L_j^{2\nu +1}}{(x^2-(K_j-1)^2)^{\frac\nu2+\frac14}}
\\ 
&&\times
G_{2\nu}(\gamma)
e\Big( 
-\frac{\eta}{2}\sqrt{x^2 - (K_j-1)^2}
- 
\frac{\eta(K_j-1)}{2\pi}
\sin^{-1}\Big( \frac{K_j-1}{x}   \Big)
\Big),
\end{eqnarray*} 
\begin{eqnarray*}\label{expansion-G-v}
G_{2\nu}(t_0)
&=&
\frac{g^{(2\nu)}(t_0)}{(2\nu)!}
+
\sum_{\ell=0}^{2\nu-1}
\frac{g^{(\ell)}(t_0)}{\ell!}
\sum_{k=1}^{2\nu-\ell}
\frac{ C_{2\nu,\ell,k} }
{\big(-\frac{\eta\pi}{2L_j^2}\sqrt{x^2- (K_j-1)^2}\big)^k}
\\
&&\times
\sum_{\mbox{\tiny$
\begin{array}{c}
3\leq n_1,n_2,\ldots,n_k\leq 2n+3
\\
n_1+n_2+\cdots+n_k=2\nu-\ell+2k
\end{array}
$}}
\frac{f^{(n_1)}(t_0)}{n_1!} 
\frac{f^{(n_2)}(t_0)}{n_2!} 
\cdots 
\frac{f^{(n_k)}(t_0)}{n_k!},
\nonumber
\end{eqnarray*}
where $\gamma = \frac{\eta L_j}{2\pi}\sin^{-1}(\frac{K_j-1}x)$ 
and $C_{2\nu,\ell,k}$ are constants. 
\end{lemma}

By Lemma \ref{Wasymp} (1), $W_{K_j,L_j}(\eta x)$ is negligible 
if $c_j\geq \frac{X}{K_j^{1-\varepsilon}L_j}$. Consequently, 
\eqref{D01} is negligible if $K_2L_2\geq X^{1+\varepsilon}$, 
\eqref{D10} is negligible if $K_1L_1\geq X^{1+\varepsilon}$, 
and \eqref{D11} is negligible if either $K_1L_1$ or 
$K_2L_2\geq X^{1+\varepsilon}$. In all other cases, we 
may restrict the $c_j$-sum to 
$c_j\leq \frac{X}{K_j^{1-\varepsilon}L_j}$. 

We will summarize the general strategy we use to bound 
$\sum\nolimits_{K_1L_1}^{K_2L_2}$ in each case. For the first 
step we use \ref{Wasymp} to rewrite the sum over $k_1$ and $k_2$ 
as an asymptotic expansion of an oscillatory integral, and we 
focus on the main term in this expansion. We then open up the 
Kloosterman sum over $z$, interchange the summation and use the 
orthogonality relation for characters to get a relation between 
the sum over $n$ and $z$. Applying Poisson summation on 
$n$, we obtain another oscillatory integral and apply a weighted 
first derivative test as in McKee, Sun, and Ye \cite{MSY} 
to shorten the sum over 
$n$. We follow a similar approach for the sum over $m$.

After applying Poisson summation for $n$ and $m$, we combine 
the two oscillatory integrals into a double integral that can 
be bounded using the following two-dimensional second derivative 
test. This is the final step in obtaining the non-trivial upper 
bound for $\sum\nolimits_{K_1L_1}^{K_2L_2}$. 

\begin{lemma}\label{2ndDer} 
{\rm (Aggarwal \cite{A} and Srinivasan \cite{Sr})} 
Suppose in a region $D\subset \mathbb{R}^2$ we have 
\begin{equation}\label{r1r2B}
\Big| \frac{\partial^2 \theta}{\partial u^2} \Big| 
\gg r_1^2, 
\ \ 
\Big| \frac{\partial^2 \theta}{\partial v^2} \Big| 
\gg r_2^2,
\ \ 
\Big| 
\frac{\partial^2 \theta}{\partial u^2}
\frac{\partial^2 \theta}{\partial v^2} 
- 
\Big(\frac{\partial^2 \theta}{\partial u\partial v}\Big)^2    
\Big|
\gg r_1^2 r_2^2,
\end{equation}
for some $r_1,r_2>0$. Define 
$$
{\rm var}(a) 
= 
\int \int_{D} 
\Big| \frac{\partial^2 a }{\partial u \partial v}\Big| 
du dv.
$$
Then 
$$
\int\int_{D} 
a(u,v)
e(\theta(u,v))dudv 
\ll 
\frac{{\rm var}(a)}{r_1 r_2}.
$$
\end{lemma}

\setcounter{equation}{0}
\section{The $D_{01}$ and $D_{10}$ terms} 

By Lemma \ref{Wasymp}, we have 
$$
D_{01}^{\eta}
=K_1K_2
\sum_{2 | k_1}
g_0\Big(\frac{k_1-K_1}{L_1}\Big)
\sum_n
\Big| \phi\Big(\frac{n}{X}\Big)\Big|^2
\sum_{c_2\leq\frac{X}{K_2^{1-\varepsilon}L_2}}
\frac{S(n,n,c_2)}{c_2}
W_{K_2,L_2}(\eta x).
$$
Notice $c_2$ is a positive integer and so we can assume 
$K_2L_2\leq X^{1+\varepsilon}$, otherwise 
$\frac{X}{K_2^{1-\varepsilon}L_2}<1$ and there will be no 
$D_{01}$ term. 
Focusing only on the leading term in the expansion for 
$\widetilde{W}_{2,0}(\eta x)$, and expanding the Kloosterman 
sum, we obtain a term of the form 
\begin{eqnarray*}
T_{01}^\eta 
& = & 
K_1K_2L_2
\sum_{2 | k_1} 
g_0\Big(\frac{k_1-K_1}{L_1}\Big) 
\sum_{n \geq 1} 
\Big|\phi\Big(\frac{n}{X}\Big)\Big|^2
\sum_{c_2 \leq  \frac{X}{K_2^{1-\varepsilon}L_2 }  }
\frac{1}{c_2}
\\
&&\times
\sideset{}{^{*}}\sum_{z \text{ mod}\,c_2}  
e\Big(\frac{nz+n\bar{z}}{c_2}\Big)
h_2^{\eta}(n,n,c_2)
e(\varphi_2^{\eta} (n,n,c_2)),
\end{eqnarray*} 
where 
\begin{eqnarray}   
\varphi_j^\eta(m, n, c) 
&=& 
-
\frac{\eta}{2\pi} 
\sqrt{\frac{16 \pi^2 mn}{c^2}-(K_j-1)^2}
-
\frac{\eta(K_j-1)}{2\pi}
\sin^{-1}\Big( \frac{(K_j-1)c}{4\pi \sqrt{mn}} \Big), 
\label{phiDef}
\\
h_j^\eta(m,n,c) 
&=& 
\hat{g}_0\Big( 
\frac{\eta L}{2\pi} 
\sin^{-1}\Big( \frac{(K_j-1)c}{4\pi \sqrt{mn}}   \Big) 
\Big) 
\Big( 
\frac{16\pi^2 mn}{c^2}-(K_j-1)^2 
\Big)^{-\frac{1}{4}}.
\label{hDef}
\end{eqnarray} 

Rewriting $n$ as $nc_2+r$, 
\begin{eqnarray}\label{T01}
T_{01}^\eta 
& = & 
K_1K_2L_2
\sum_{2 | k_1} 
g_0\Big(\frac{k_1-K_1}{L_1}\Big) 
\sum_{c_2 \leq  \frac{X}{K_2^{1-\varepsilon}L_2 }  }
\frac{1}{c_2} 
\sideset{}{^{*}}\sum_{z \text{ mod}\,c_2}
\sum_{r \text{ mod}\,c_2}  
e\Big(\frac{rz+r\bar{z}}{c_2}\Big) 
\\
&&\times 
\sum_{n \geq 1} 
\Big|\phi\Big(\frac{nc_2+r}{X}\Big)\Big|^2 
h_2^{\eta}(nc_2+r,nc_2+r,c_2)
e(\varphi_2^{\eta} (nc_2+r,nc_2+r,c_2)).
\nonumber
\end{eqnarray} 
Applying the Poisson summation to the $n$-sum in 
\eqref{T01}, we get 
\begin{eqnarray}\label{nSum}
&&
\sum_{n\in\mathbb Z}
\int_{\mathbb{R}} 
\Big|\phi\Big(\frac{yc_2+r}{X}\Big)\Big|^2 
h_2^{\eta}(yc_2+r,yc_2+r,c_2)
e(\varphi_2^{\eta} (yc_2+r,yc_2+r,c_2))
e(-yn)dy
\\
&=&
\frac1{c_2}
\sum_{n\in\mathbb Z} 
e\Big(\frac{rn}{c_2}\Big) 
\int_{\mathbb{R}} 
\Big|\phi\Big(\frac{t}{X}\Big)\Big|^2 
h_2^{\eta}(t,t,c_2)
e\Big(\varphi_2^{\eta} (t,t,c_2)-\frac{tn}{c_2}\Big)
dt
\nonumber
\end{eqnarray}
by changing variables to $t=yc_2+r$. Substituting \eqref{nSum} 
back to \eqref{T01}, we may evaluating the $r$-sum and get rid 
of the Kloosterman sum
\begin{eqnarray}\label{T01afterP}
T_{01}^\eta 
& = & 
K_1K_2L_2
\sum_{2 | k_1} 
g_0\Big(\frac{k_1-K_1}{L_1}\Big) 
\sum_{c_2 \leq  \frac{X}{K_2^{1-\varepsilon}L_2 }  }
\frac{1}{c_2} 
\sideset{}{^{*}}\sum_{z \text{ mod}\,c_2}
\sum_{n \equiv -z-\bar{z}\text{ mod}\,c_2} 
\\
&&\times 
\int_{\mathbb{R}} 
\Big|\phi\Big(\frac{t}{X}\Big)\Big|^2 
h_2^{\eta}(t,t,c_2)
e\Big(\varphi_2^{\eta} (t,t,c_2)-\frac{tn}{c_2}\Big)
dt.
\nonumber
\end{eqnarray}

We denote by $I$ the integral on the right hand side of 
\eqref{T01afterP} and change variables from $t$ to $tX$, 
thus making it an integral over $[1,2]$, the support of $\phi$
$$
I
=
X
\int_{1}^2 
|\phi(t)|^2 
h_2^{\eta}(Xt,Xt,c_2)
e\Big(\varphi_2^{\eta} (Xt,Xt,c_2)-\frac{tnX}{c_2}\Big)
dt.
$$
By the Riemann-Lebesgue lemma $I$ is negligible for $n$ 
outside a compact interval. 
Note that $\frac{d^l\phi}{dt^l}\ll 1$ and by 
\eqref{hDef} 
\begin{eqnarray*}
h_2^\eta(Xt,Xt,c_2) 
&=&
\hat{g}_0\Big( 
\frac{\eta L_2}{2\pi} 
\sin^{-1}\Big( \frac{(K_2-1)c_2}{4\pi Xt}\Big) 
\Big) 
\Big( 
\frac{16\pi^2 X^2t^2}{c_2^2}-(K_2-1)^2 
\Big)^{-\frac{1}{4}}
\\
&\asymp&
\Big(\frac{X^2}{c_2^2}\Big)^{-\frac{1}{4}}
\asymp 
c_2^\frac12X^{-\frac12}
=:U
\end{eqnarray*}
as $16\pi^2X^2t^2c_2^{-2}$ always dominates $(K_2-1)^2$. 
Taking derivatives,
\begin{eqnarray*}
&&
\frac{d}{dt}
\Big( 
\frac{16\pi^2 X^2t^2}{c_2^2}-(K_2-1)^2
\Big)^{-\frac{1}{4}} 
\asymp U,
\\
&&
\frac{d}{dt}
\hat{g}_0\Big( 
\frac{\eta L_2}{2\pi} 
\sin^{-1}\Big( \frac{(K_2-1)c_2}{4\pi Xt}\Big)
\Big) 
\\
& = & 
\hat{g}'_0\Big( 
\frac{\eta L_2}{2\pi} 
\sin^{-1}\Big( \frac{(K_2-1)c_2}{4\pi Xt} \Big) 
\Big)
\frac{\eta L_2}{2\pi\sqrt{1-(\frac{(K_2-1)c_2}{4\pi Xt})^2}}
\frac{-(K_2-1)c_2}{4\pi Xt^2}
\asymp
\frac{K_2L_2c_2}{X}
\ll 
K_2^{\varepsilon}.
\end{eqnarray*} 
Subsequent derivatives yield $\ll K_2^{\varepsilon}$ too.
Putting them together we have 
\begin{equation}\label{phihDer}
\frac{d^l}{dt^l}
|\phi(t)|^2
h_2^{\eta}(Xt, Xt, c_2)
\ll \frac{U}{N^l},
\ l\geq 0,
\end{equation}
for $N=K_2^{-\varepsilon}$. Recall the definition of 
$\varphi_2^{\eta}(Xt, Xt, c_2)$ in \eqref{phiDef}. 
Computing derivatives we get 
\begin{eqnarray}
\frac{d}{dt}
\varphi_2^\eta(Xt, Xt, c_2)
-\frac{tnX}{c_2} 
& = & 
-\frac{\eta}{4\pi} 
\frac{\frac{16\pi^2X^2 2t}{c_2^2}}
{\sqrt{\frac{16 \pi^2 X^2 2t}{c^2}-(K_2-1)^2}}
-\frac{\eta(K_2-1)}{2\pi}
\frac{-\frac{(K_2-1)c_2}{4\pi Xt^2}}
{\sqrt{1-\frac{(K_2-1)^2c_2^2}{(4\pi)^2X^2t^2}}}
-\frac{nX}{c_2}
\nonumber
\\
& = &
-\frac{\eta}{2\pi t}
\sqrt{\frac{16\pi^2X^2t^2}{c_2^2}-(K_2-1)^2}
-\frac{nX}{c_2}
\nonumber
\\
& = & 
-\frac{4\eta\pi Xt}{2\pi tc_2}
\Big(1+O\Big(\frac{c_2^2K_2^2}{X^2}\Big)\Big)
-\frac{n X}{c_2}
=
-(2\eta +n)
\frac{X}{c_2}
\Big(1+O\Big(\frac{c_2^2K_2^2}{X^2}\Big)\Big), 
\label{phiDer}
\end{eqnarray}
since 
$$
\frac{c_2^2K_2^2}{X^2}
\leq 
\frac{K_2^2}{K_2^{2-\varepsilon}L_2^2}
=
\frac{K_2^{\varepsilon}}{L_2^2} 
$$ 
is small. Now from \eqref{phiDer} using \eqref{phihDer}, 
if $n\neq -2\eta$, we have 
$$
\Big|
\frac{d}{dt}
\varphi_2^{\eta}(Xt,Xt, c_2)
-\frac{tnX}{c_2}
\Big|
\gg 
\frac{X}{c_2}
\geq K_2^{1-\varepsilon}L_2.
$$
Subsequent derivatives are all $\asymp Xc_2^{-1}$. Thus we may 
take $T=Xc_2^{-1}$ and $M=1$. By the first derivative test 
(cf. \cite{MSY}), the integral 
$$
I
\ll 
U
\Big(\frac{c_2K_2^\varepsilon}X\Big)^{n_0+1}
\leq
\frac{c_2^\frac12}{X^\frac12(K_2^{1-2\varepsilon}L_2)^{n_0+1}} 
$$ 
is negligible for $n_0$ sufficiently large, when $n\neq -2\eta$.

Excluding these negligible terms we get 
\begin{eqnarray*}
T_{01}^{\eta}
& = & 
K_1K_2L_2
\sum_{2|k_1}
g_0\Big(\frac{k_1-K_1}{L_1}\Big)
\sum_{c_2\leq\frac{X}{K_2^{1-\varepsilon}L_2}}
\frac{1}{c_2}
\sideset{}{^{*}}\sum_{
\genfrac{}{}{0pt}{1}{z\,{\rm mod}\,c_2}
{z+\bar{z}\equiv 2\eta\,{\rm mod}\,c_2}
}
\\
&&\times 
\int_{\mathbb{R}}
\Big|\phi\Big(\frac{t}{X}\Big)\Big|^2
h_2^{\eta}(t,t,c_2)
e\Big(\varphi_2^{\eta}(t, t, c_2)+\frac{2\eta t}{c_2}\Big)
dt
+ O(K_2^{-A}).
\end{eqnarray*}
The congruence $z+\bar{z}\equiv 2\eta \text{ mod } c_2$ is 
equivalent to $z^2+1\equiv 2\eta z \text{ mod } c_2$ or 
equivalently to $(z-\eta)^2\equiv 0 \text{ mod } c_2. $
For $c_2=p^c$ with $c\geq 1$, this means $p^c | (z-\eta)^2$, 
i.e.,
$p^{\frac{c}{2}} | z-\eta$ if c is even, and 
$p^{\frac{c+1}{2}} | z-\eta$ if $c$ is odd. Consequently 
the number of solutions of 
$z+\bar{z}\equiv 2\eta \text{ mod } c_2$ for 
$z \text{ mod } c_2,\ (z,c_2)=1$, is $\leq \sqrt{c_2}$.
Now we compute the integral $I$ for $n=-2 \eta$.

By \eqref{phiDer}
\begin{eqnarray*}
&&
\frac{d}{dt}
\varphi_2^{\eta}(Xt,Xt, c_2)
+
\frac{2\eta tX}{c_2} 
= 
-\frac{\eta}{2\pi t}
\sqrt{\frac{16\pi^2 X^2t^2}{c_2^2}-(K_2-1)^2}
+\frac{2\eta X}{c_2}
\\
& = & 
\frac{
\frac{4X^2}{c_2^2}
-\frac{1}{4\pi^2t^2}
\Big(\frac{16\pi^2 X^2t^2}{c_2^2}-(K_2-1)^2\Big)
}{
\frac{\eta}{2 \pi t}
\sqrt{\frac{16\pi^2 X^2t^2}{c_2^2}-(K_2-1)^2}
+\frac{2\eta X}{c_2}
}
= 
\frac{
\frac{(K_2-1)^2}{4\pi^2t^2}
}{
\frac{\eta}{2 \pi t}
\sqrt{\frac{16\pi^2 X^2t^2}{c_2^2}-(K_2-1)^2}
+\frac{2\eta X}{c_2}}
\asymp 
\frac{c_2K_2^2}{X}.
\nonumber
\end{eqnarray*}
Any subsequent differentiation yields a factor of $\asymp 1$. 
By the first derivative test again, 
$$
I
\ll 
U
\Big(\frac X{c_2K_2^{2-\varepsilon}}\Big)^{n_0+1}
$$ 
which is negligible for $n_0$ sufficiently large if 
$\frac X{c_2K_2^{2-\varepsilon}}\leq K_2^{-\varepsilon}$, 
i.e., if 
$$
\frac{X}{K_2^{2-\varepsilon}}
\leq c_2\leq 
\frac{X}{K_2^{1-\varepsilon}L_2}. 
$$
Therefore we can shorten the $c_2$-sum and reduce 
$T_{01}^{\eta}$ to
\begin{eqnarray*}
T_{01}^{\eta} 
& = & 
K_1K_2L_2
\sum_{2|k_1}
g_0\Big(\frac{k_1-K_1}{L_1}\Big)
\sum_{c_2\leq\frac{X}{K_2^{2-\varepsilon}}}
\frac{1}{c_2}
\sideset{}{^{*}}\sum_{
\genfrac{}{}{0pt}{1}{z\,{\rm mod}\,c_2}
{z+\bar{z}\equiv 2\eta\,{\rm mod}\,c_2}
}
\\
&&\times 
\int_{\mathbb{R}}
\Big|\phi\Big(\frac{t}{X}\Big)\Big|^2
h_2^{\eta}(t,t,c_2)
e\Big(
\varphi_2^{\eta}(t, t, c_2)
+\frac{2\eta t}{c_2}
\Big)dt
+ O(K_2^{-A}),
\end{eqnarray*}
if $K_2\leq X^{\frac12+\varepsilon}$, while $T_{01}^{\eta}$ 
is negligible if $K_2\geq T^{\frac{1}{2}+\varepsilon}$. 

By trivial estimation, we have $I\ll c_2^\frac12X^\frac12$. 
Consequently, 
\begin{eqnarray}\label{T01B}
T_{01}^{\eta} 
& \ll & 
K_1L_1K_2L_2X^{\frac{1}{2}}
\sum_{c_2\leq \frac{X}{K_2^{2-\varepsilon}}}
\frac{1}{c_2}c_2^{\frac{1}{2}}c_2^{\frac{1}{2}} 
+O(K_2^{-A})
\\
& \ll & 
\frac{K_1L_1L_2X^{\frac{3}{2}}}{K_2^{1-\varepsilon}}
\ \ {\rm if}\ K_2\leq X^{\frac{1}{2}+\varepsilon} 
\nonumber
\\
& \ll & 
K_2^{-A}
\ \ {\rm if}\ K_2\geq X^{\frac{1}{2}+\varepsilon}.
\nonumber
\end{eqnarray}
Similarly,
\begin{eqnarray}\label{T10B}
T_{10}^{\eta} 
& \ll & 
\frac{K_2L_1L_2X^{\frac{3}{2}}}{K_1^{1-\varepsilon}}
\ \ {\rm if}\ K_1\leq X^{\frac{1}{2}+\varepsilon} 
\\
& \ll & 
K_1^{-A}
\ \ {\rm if}\ K_1\geq X^{\frac{1}{2}+\varepsilon}.
\nonumber
\end{eqnarray}

These give upper bounds for the main terms in the 
expansions of $D_{01}$ and $D_{10}$ for $0<\beta\leq1$. 
Bounds in \eqref{K1K2>} and \eqref{K1>K2<} then follow. 

\setcounter{equation}{0}
\section{Poisson summation for $D_{11}$} 

To use Poisson summation, we rewrite for $0<\beta\leq1$ 
\begin{eqnarray*}
T_{11}^{\eta_1\eta_2} 
& = & 
K_1L_1K_2L_2
\sum_n
\sum_m 
e(\alpha n^{\beta}-\alpha m^{\beta})
\phi(\frac{n}{X})
\bar{\phi}(\frac{m}{X})
\sum_{d\leq\min(\frac{X}{K_1^{1-\varepsilon}L_1}, 
\frac{X}{K_2^{1-\varepsilon}L_2})}
\\
&&\times
\sum_{
\genfrac{}{}{0pt}{1}{c_j\leq\frac X{dK_j^{1-\varepsilon}L_j}}
{(c_1,c_2)=1}
}
\frac{1}{c_1 c_2 d^2}
\sideset{}{^{*}}\sum_{z_1\,{\rm mod}\,dc_1}
e\Big(\frac{mz_1+n\bar{z}_1}{dc_1}\Big)
\sideset{}{^{*}}\sum_{z_2\,{\rm mod}\,dc_2}
e\Big(\frac{mz_2+n\bar{z}_2}{dc_2}\Big)
\nonumber
\\
&&\times
h_1^{\eta_1}(m, n, dc_1)
h_2^{\eta_2}(m, n, dc_2)
e(\varphi_1^{\eta_1}(m,n,dc_1)+\varphi_2^{\eta_2}(m,n,dc_2)).
\nonumber
\end{eqnarray*}
Rewriting $n$ as 
$nc_1c_2d+r$, we have
\begin{eqnarray}\label{T11Change}
T_{11}^{\eta_1\eta_2} 
& = & 
K_1L_1K_2L_2
\sum_m 
e(-\alpha m^{\beta}) 
\bar{\phi}(\frac{m}{X})
\sum_{d\leq\min(\frac{X}{K_1^{1-\varepsilon}L_1}, 
\frac{X}{K_2^{1-\varepsilon}L_2})}
\sum_{
\genfrac{}{}{0pt}{1}{c_j\leq\frac X{dK_j^{1-\varepsilon}L_j}}
{(c_1,c_2)=1}
}
\frac{1}{c_1c_2d^2}
\\
&&\times
\sideset{}{^{*}}\sum_{z_1\,{\rm mod}\,dc_1}
e\Big(\frac{mz_1}{dc_1}\Big)
\sideset{}{^{*}}\sum_{z_2\,{\rm mod}\,dc_2}
e\Big(\frac{mz_2}{dc_2}\Big)
\sum_{r\,{\rm mod}\,dc_1c_2}
e\Big(\frac{r\bar{z}_1}{dc_1}+\frac{r\bar{z}_2}{dc_2} \Big)
\nonumber
\\
&&\times
\sum_n
\phi\Big(\frac{nc_1c_2d+r}{X}\Big)
h_1^{\eta_1}(m, nc_1c_2d+r, c_1d)
h_2^{\eta_2}(m, nc_1c_2d+r, c_2d)
\nonumber
\\
&&\times
e(
\alpha(nc_1c_2d+r)^{\beta} 
+
\varphi_1^{\eta_1}(m, nc_1c_2d+r, c_1d)
+
\varphi_2^{\eta_2}(m, nc_1c_2d+r, c_2d)
).
\nonumber
\end{eqnarray}
By Poisson summation on the $n$-sum, we get
\begin{eqnarray*}
&&
\sum_n
\int_{\mathbb{R}}
\phi\Big(\frac{yc_1c_2d+r}{X}\Big)
h_1^{\eta_1}(m, yc_1c_2d+r, c_1d)
h_2^{\eta_2}(m, yc_1c_2d+r, c_2d)
\\
&&\times
e(
\alpha(yc_1c_2d+r)^{\beta}
+\varphi_1^{\eta_1}(m, yc_1c_2d+r, c_1d)
+\varphi_2^{\eta_2}(m, yc_1c_2d+r, c_2d)
-ny
)dy
\\
& = & 
\sum_n
\frac{1}{c_1c_2d}
e\Big(\frac{nr}{c_1c_2d}\Big)
\int_{\mathbb{R}}
\phi\Big(\frac{v}{X}\Big)
h_1^{\eta_1}(m, v, c_1d)
h_2^{\eta_2}(m, v, c_2d)
\\
&&\times
e\Big(
\alpha v^{\beta}
+\varphi_1^{\eta_1}(m, v, c_1d)
+\varphi_2^{\eta_2}(m, v, c_2d)
-\frac{nv}{c_1c_2d}
\Big)dv.
\end{eqnarray*}

The $r$-sum in \eqref{T11Change} becomes
$$
\sum_{r\,{\rm mod}\,c_1c_2d}
e\Big(\frac{r(n+c_2\bar{z}_1+c_1\bar{z}_2)}{c_1c_2d}\Big)= 
\begin{cases}
    c_1c_2d & \text{if } 
n\equiv -c_2\bar{z}_1-c_1\bar{z}_2 \text{ mod } c_1c_2d\\
    0              & \text{otherwise}.
\end{cases}
$$
Consequently    
\begin{eqnarray*}
T_{11}^{\eta_1\eta_2} 
& = & 
K_1L_1K_2L_2
\sum_{m}
e(-\alpha m^{\beta})
\bar{\phi}\Big(\frac mX\Big)
\sum_{d\leq\min(\frac{X}{K_1^{1-\varepsilon}L_1}, 
\frac{X}{K_2^{1-\varepsilon}L_2})}
\frac{1}{d^2}
\sum_{
\genfrac{}{}{0pt}{1}{c_j\leq\frac X{dK_j^{1-\varepsilon}L_j}}
{(c_1,c_2)=1}
}
\frac{1}{c_1c_2}
\\
&&\times
\sideset{}{^{*}}\sum_{z_1\,{\rm mod}\,dc_1}
e\Big(\frac{mz_1}{dc_1}\Big)
\sideset{}{^{*}}\sum_{z_2\,{\rm mod}\,dc_2}
e\Big(\frac{mz_2}{dc_2}\Big)
\sum_{n\equiv -c_2\bar{z}_1-c_1\bar{z}_2\,({\rm mod}\,c_1c_2d)}
\int_{\mathbb{R}}
\phi\Big(\frac{v}{X}\Big)
h_1^{\eta_1}(m, v, c_1d)
\\
&&\times
h_2^{\eta_2}(m, v, c_2d)
e\Big(
\alpha v^{\beta}
+\varphi_1^{\eta_1}(m, v, c_1d)
+\varphi_2^{\eta_2}(m, v, c_2d)
-\frac{nv}{c_1c_2d}
\Big)dv.
\end{eqnarray*}
Applying the same to the $m$-sum, we have
\begin{eqnarray}\label{T11afterP}
T_{11}^{\eta_1\eta_2} 
& = & 
K_1L_1K_2L_2
\sum_{d\leq\min(\frac{X}{K_1^{1-\varepsilon}L_1}, 
\frac{X}{K_2^{1-\varepsilon}L_2})}
\frac{1}{d^2}
\sum_{
\genfrac{}{}{0pt}{1}{c_j\leq\frac X{dK_j^{1-\varepsilon}L_j}}
{(c_1,c_2)=1}
}
\frac{1}{c_1c_2}
\sideset{}{^{*}}\sum_{z_1\,{\rm mod}\,dc_1}
\ 
\sideset{}{^{*}}\sum_{z_2\,{\rm mod}\,dc_2}
\\
& \times & 
\sum_{m\equiv -c_2z_1-c_1z_2\,({\rm mod}\,c_1c_2d)}
\ 
\sum_{n\equiv -c_2\bar{z}_1-c_1\bar{z}_2\,({\rm mod}\,c_1c_2d)}
\int_{\mathbb{R}}
\int_{\mathbb{R}}
a(u,v)
e(\theta^{\eta_1\eta_2}(u,v))
dudv,
\nonumber
\end{eqnarray}
where 
\begin{eqnarray*}
a(u,v)
&=&
\phi\Big(\frac{v}{X}\Big)
\bar{\phi}\Big(\frac{u}{X}\Big)
h_1^{\eta_1}(u,v,c_1d)
h_2^{\eta_2}(u,v,c_2d), 
\\
\theta^{\eta_1\eta_2}(u,v)
&=&
\alpha v^{\beta}
-\alpha u^{\beta}
+\varphi_1^{\eta_1}(u,v,c_1d)
+\varphi_2^{\eta_2}(u,v,c_2d)
-\frac{mu}{c_1c_2d}
-\frac{nv}{c_1c_2d}
=:
\theta(u,v). 
\end{eqnarray*}

By \eqref{hDef}, $a(u,v)\ll dc_1^\frac12c_2^\frac12X^{-1}=:U$, 
and each $\frac{\partial}{\partial u}$ produces a factor 
$\ll K_1^{\varepsilon}K_2^{\varepsilon}X^{-1}=\frac{1}{N}$ with 
$N=XK_1^{-\varepsilon}K_2^{-\varepsilon}$. Computing derivatives 
in the case of $0<\beta<1$ we get 
\begin{equation}\label{dutheta}
\frac{\partial\theta^{\eta_1\eta_2}}{\partial u}
=
-\alpha\beta u^{\beta-1}
-\frac{\eta_1}{4\pi u}
\sqrt{R_1}
-\frac{\eta_2}{4\pi u}
\sqrt{R_2}
-\frac{m}{c_1c_2d},
\end{equation}
where $R_j=\frac{16\pi^2 uv}{c_j^2d^2}-(K_j-1)^2$. 
Suppose that the absolute value of the sum of the first three 
terms on the right hand side of \eqref{dutheta} is bounded by 
$\frac{\tau-\delta}{c_1c_2d}$ for some $0<\delta<\tau$. 
Then for 
$|m|\geq \tau$, we have  
$|\frac{\partial\theta}{\partial u}|\geq \frac\delta{c_1c_2d}$. 
Moreover for $r\geq 2$, 
$\frac{\partial^r\theta}{\partial u^r}\ll \frac{T}{M^r}$,
where $T=\frac{X\tau}{c_1c_2d}$ and $M=X$. 
Denote by $J$ the double integral in \eqref{T11afterP}. Then 
by the first derivative test \cite[Theorem 1.1]{MSY}, 
$J$ is negligible for 
$|m|\geq\tau$. Therefore, 
\begin{eqnarray}\label{T11J}
T_{11}^{\eta_1\eta_2} 
& = & 
K_1L_1K_2L_2
\sum_{d\leq\min(\frac{X}{K_1^{1-\varepsilon}L_1}, 
\frac{X}{K_2^{1-\varepsilon}L_2})}
\frac{1}{d^2} 
\sum_{
\genfrac{}{}{0pt}{1}{c_j\leq\frac X{dK_j^{1-\varepsilon}L_j}}
{(c_1,c_2)=1}
}
\frac{1}{c_1c_2}
\sideset{}{^{*}}\sum_{z_1\,{\rm mod}\,dc_1}
\ 
\sideset{}{^{*}}\sum_{z_2\,{\rm mod}\,dc_2}
\\
& \times & 
\sum_{
\genfrac{}{}{0pt}{1}{|m|\leq \tau}
{m\equiv-c_2z_1-c_1z_2\,({\rm mod}\,c_1c_2d)}
} 
\ 
\sum_{
\genfrac{}{}{0pt}{1}{|n|\leq \tau}
{n\equiv-c_2\bar{z}_1-c_1\bar{z}_2\,({\rm mod}\,c_1c_2d)}
}
J
+O(X^{-A}).
\nonumber
\end{eqnarray}

\begin{lemma}\label{congr}
The $z_1$, $z_2$, $m$, $n$-sums in \eqref{T11J} have at most 
$d(2\tau+1)([\frac\tau{c_1c_2d}]+1)$ terms. 
\end{lemma}

{\it Proof.} 
Reducing 
\begin{equation}\label{mCongr}
m\equiv-c_2z_1-c_1z_2\,({\rm mod}\,c_1c_2d)
\end{equation}
to congruences mod $c_1$ and mod $c_2$, we see that 
for each $m$, 
$z_j$ is uniquely determined ${\rm mod}\,c_j$, $j=1,2$.
But modulo $c_jd$, there are $d$ such $z_j$'s: $z_j+c_jk_j$ with 
$0\leq k_j<d$. Then \eqref{mCongr} becomes 
$$
k_1+k_2\equiv
-\frac{m+c_1z_2+c_2z_1}{c_1c_2}
\,({\rm mod}\,d).
$$
Consequently, given $m$, there are at most $d$ such 
$z_1\,{\rm mod}\,dc_1$. Given $m$ and $z_1$, there is a unique 
$z_2\,{\rm mod}\,dc_2$. Given $m$, $z_1$, and $z_2$, there 
is a unique $n\,{\rm mod}\,dc_1c_2$. The lemma then follows 
because there are at most $2\tau+1$ $m$'s. 
\hfill
$\square$

We want to apply the second derivative test to get an upper 
bound for the double integral $J$. We have the following 
second derivatives in the case of $0<\beta<1$:
\begin{eqnarray}
\frac{\partial^2\theta^{\eta_1\eta_2}}{\partial u^2} 
& = & 
-\alpha\beta(\beta-1)u^{\beta-2}
+
\frac{\eta_1\sqrt{R_1}}{4\pi u^2}
+
\frac{\eta_2\sqrt{R_2}}{4\pi u^2}
-
\frac{2\pi\eta_1v}{uc_1^2d^2\sqrt{R_1}}
-
\frac{2\pi\eta_2v}{uc_2^2d^2\sqrt{R_2}}
=:
\sum_{i=1}^5U_i,
\label{uu}
\\
\frac{\partial^2\theta^{\eta_1\eta_2}}{\partial v^2} 
& = & 
\alpha\beta(\beta-1)v^{\beta-2}
+
\frac{\eta_1\sqrt{R_1}}{4\pi v^2}
+
\frac{\eta_2\sqrt{R_2}}{4\pi v^2}
-
\frac{2\pi\eta_1u}{vc_1^2d^2\sqrt{R_1}}
-
\frac{2\pi\eta_2u}{vc_2^2d^2\sqrt{R_2}}
=:
\sum_{i=1}^5V_i,
\label{vv}
\\
\frac{\partial^2\theta^{\eta_1\eta_2}}{\partial u\partial v} 
& = & 
-
\frac{2\pi\eta_1}{c_1^2d^2\sqrt{R_1}}
-
\frac{2\pi\eta_2}{c_2^2d^2\sqrt{R_2}}.
\nonumber
\end{eqnarray}
Consequently, 
\begin{eqnarray}
\Big(\frac{\partial^2\theta}{\partial u\partial v}\Big)^2 
& = & 
(U_4+U_5)(V_4+V_5),
\nonumber
\\
\frac{\partial^2\theta}{\partial u^2}
\frac{\partial^2\theta}{\partial v^2}
-
\Big(\frac{\partial^2\theta}{\partial u\partial v}\Big)^2 
& = & 
U_1V_1
+
U_1\sum_{i=2}^5V_i
+
V_1\sum_{i=2}^5U_i
\label{-uvSq}
\\
&&+
(U_2+U_3)\sum_{i=2}^5V_i
+
(U_4+U_5)(V_2+V_3).
\label{UV2-5}
\end{eqnarray}
Note that \eqref{UV2-5} equals
\begin{eqnarray}
&&
\frac{(\eta_1\sqrt{R_1}+\eta_2\sqrt{R_2})^2}{16\pi^2u^2v^2}
-
\frac{\eta_1\sqrt{R_1}+\eta_2\sqrt{R_2}}{uvd^2}
\Big(
\frac{\eta_1}{c_1^2\sqrt{R_1}}
+
\frac{\eta_2}{c_2^2\sqrt{R_2}}
\Big)
\label{R1-R2}
\\
&=&
-
\frac{\eta_1\sqrt{R_1}+\eta_2\sqrt{R_2}}{16\pi^2u^2v^2}
\Big(
\frac{\eta_1(K_1-1)^2}{\sqrt{R_1}}
+
\frac{\eta_2(K_2-1)^2}{\sqrt{R_2}}
\Big)
\label{R1+R2}
\end{eqnarray}
by
$$
\frac{\sqrt{R_j}}{16\pi^2u^2v^2}
-
\frac1{uvc_j^2d^2\sqrt{R_j}}
=
\frac{\frac{16\pi^2uv}{c_j^2d^2}-(K_j-1)^2}
{16\pi^2u^2v^2\sqrt{R_j}}
-
\frac1{uvc_j^2d^2\sqrt{R_j}}
=
-\frac{(K_j-1)^2}{16\pi^2u^2v^2\sqrt{R_j}}.
$$

\setcounter{equation}{0}
\section{Case of $\eta_1=\eta_2$ for $D_{11}$ when $0<\beta<1$} 

Assume $\eta_1=\eta_2$ in this section. Then there is no 
cancellation in the middle two terms of \eqref{dutheta}, and 
hence their sum is equal to
$$
-\frac{\eta_1\sqrt{R_1}}{4\pi u}
-\frac{\eta_2\sqrt{R_2}}{4\pi u}
=
-\frac{\eta_1}d
\sqrt{\frac vu}
\Big(\frac1{c_1}+\frac1{c_2}\Big)
+O\Big(\frac{c_1dK_1^2}{X^2}\Big)
+O\Big(\frac{c_2dK_2^2}{X^2}\Big),
$$
where the first term on the right hand side dominates. 
Thus,
\begin{equation}\label{Mid2Terms}
\frac{c_1+c_2}{\sqrt2dc_1c_2}
\leq
\Big|
-\frac{\eta_1\sqrt{R_1}}{4\pi u}
-\frac{\eta_2\sqrt{R_2}}{4\pi u}
\Big|
\leq
\frac{\sqrt2(c_1+c_2)}{dc_1c_2}.
\end{equation}
We will assume $K_1=K_2$ and 
\begin{equation}\label{KLlarge}
K_1^2L_1L_2\geq X^{1+\beta+\varepsilon}.
\end{equation}
Then $K_1L_1,K_1L_2\geq X^{\beta+\frac\varepsilon2}$ 
because $K_1L_1,K_1L_2\leq X^{1+\frac\varepsilon2}$. 
Consequently, 
$|-\alpha\beta u^{\beta-1}|\leq\alpha\beta X^{\beta-1}$ 
is a power smaller than 
$\frac1d(\frac1{c_1}+\frac1{c_2})$. 
Thus, the absolute value of the sum of the first three 
terms on the right hand side of \eqref{dutheta} is 
$$
\leq 
\alpha\beta X^{\beta-1}
+\frac{\sqrt{2}(c_1+c_2)}{c_1c_2d}
\leq
\frac{1.42(c_1+c_2)}{c_1c_2d}, 
$$
and we may take $\tau=1.5(c_1+c_2)$ and $\delta=0.08(c_1+c_2)$ 
and apply Lemma \ref{congr} to get \eqref{T11J}. 

Now we bound $J$ in \eqref{T11J}. Because $\eta_1=\eta_2$, 
there are no cancellations in \eqref{R1+R2}. Consequently, 
\eqref{UV2-5} is
\begin{equation}\label{UV2-5B}
\asymp 
\frac1{X^4}
\Big(\frac X{c_1d}+\frac X{c_2d}\Big)
\Big(\frac{c_1dK_1^2}X+\frac{c_2dK_2^2}X\Big)
=
\frac{K_1^2(c_1+c_2)}{X^4}
\Big(\frac1{c_1}+\frac1{c_2}\Big).
\end{equation}
By 
\begin{equation}\label{Rj-/Rj}
\frac{\sqrt{R_j}}{4\pi u^2}
-\frac{2\pi v}{uc_j^2d^2\sqrt{R_j}}
=
\frac{\frac{16\pi^2uv}{c_j^2d^2}-(K_j-1)^2}
{4\pi u^2\sqrt{R_j}}
-\frac{2\pi v}{uc_j^2d^2\sqrt{R_j}}
=
\frac1{2c_jdu}
\sqrt{\frac vu}
\Big(1+
O\Big(\frac{c^2_jd^2K_j^2}{X^2}\Big)
\Big),
\end{equation}
we know that 
$$
\sum_{i=2}^5
U_i
\asymp
\frac1{dX}
\Big(\frac1{c_1}+\frac1{c_2}\Big)
\geq
\frac{K_1^{1-\varepsilon}L_1+K_1^{1-\varepsilon}L_2}{X^2}
\gg
X^{\beta-2+\varepsilon}
\gg
U_1
$$
by \eqref{KLlarge}. Likewise, 
$\sum_{i=2}^5V_i$ dominates $V_1$. Therefore, 
\eqref{uu} and \eqref{vv} are both 
$\asymp\frac1{dX}(\frac1{c_1}+\frac1{c_2})$. By the same reason,
\begin{equation}\label{V1U}
U_1\sum_{i=1}^5V_i
\asymp
V_1\sum_{i=2}^5U_i
\asymp
\frac{X^{\beta-3}}d
\Big(\frac1{c_1}+\frac1{c_2}\Big).
\end{equation}
Since \eqref{UV2-5} dominates \eqref{V1U}, the left hand 
side of \eqref{-uvSq} is $\asymp$ \eqref{UV2-5B}. We observe 
that 
\begin{eqnarray}
\frac1{d^2X^2}
\Big(\frac1{c_1}+\frac1{c_2}\Big)^2
&\geq&
\frac{K_1^{1-\varepsilon}L_1+K_1^{1-\varepsilon}L_2}{dX^3}
\Big(\frac1{c_1}+\frac1{c_2}\Big),
\label{r1^4}
\\
\frac{K_1^2(c_1+c_2)}{X^4}
\Big(\frac1{c_1}+\frac1{c_2}\Big)
&\leq&
\frac1{dX^3}
\Big(
\frac{K_1^{1+\varepsilon}}{L_1}
+
\frac{K_1^{1+\varepsilon}}{L_2}
\Big)
\Big(\frac1{c_1}+\frac1{c_2}\Big).
\label{r1r2^2}
\end{eqnarray}
Since \eqref{r1^4} is a power larger than \eqref{r1r2^2}, 
we may choose 
$$
r_1
=
r_2
=
\frac{K_1^\frac12(c_1+c_2)^\frac14}{X}
\Big(\frac{1}{c_1}+\frac{1}{c_2}\Big)^{\frac{1}{4}}
=
\frac{K_1^\frac12(c_1+c_2)^\frac12}{(c_1c_2)^\frac14X}
$$
as in \eqref{r1r2B}.

Note that $a(u,v)\ll d(c_1c_2)^\frac12X^{-1}$ and each 
differentiation produces a factor 
$\ll K_1^{\varepsilon}X^{-1}$. 
Consequently,
\begin{equation}\label{varB}
\frac{\partial^2a}{\partial u\partial v} 
\ll 
\frac{K_1^{\varepsilon}d\sqrt{c_1c_2}}{X^3},
\ \ \ 
{\rm var}(a) 
\ll 
\frac{K_1^{\varepsilon}d\sqrt{c_1c_2}}{X}.
\end{equation}
By Lemma \ref{2ndDer} for $\eta_1=\eta_2$ we have 
$$
J 
\ll 
\frac{d\sqrt{c_1c_2}K_1^\varepsilon}X
\frac{\sqrt{c_1c_2}X^2}{(c_1+c_2)K_1}
=
\frac{c_1c_2dX}{(c_1+c_2)K_1^{1-\varepsilon}},
$$
and hence by \eqref{T11J}
\begin{eqnarray*}
T_{11}^{\eta_1\eta_2} 
& \ll & 
K_1^{1+\varepsilon}L_1L_2
X
\sum_{d\leq\min(\frac{X}{K_1^{1-\varepsilon}L_1}, 
\frac{X}{K_1^{1-\varepsilon}L_2})}
\frac{1}{d}
\sum_{\genfrac{}{}{0pt}{1}
{c_j\leq\frac{X}{dK_j^{1-\varepsilon}L_j}}{(c_1,c_2)=1}}
\ 
\sideset{}{^{*}}\sum_{z_1\,{\rm mod}\,dc_1}
\ 
\sideset{}{^{*}}\sum_{\text{  }z_2\,{\rm mod}\,dc_2}
\\
&&\times 
\sum_{\genfrac{}{}{0pt}{1}{|m|\leq 1.5(c_1+c_2)}
{m\equiv -c_1z_1-c_2z_2\,({\rm mod}\,c_1c_2d)}}
\sum_{\genfrac{}{}{0pt}{1}{|n|\leq 1.5(c_1+c_2)}
{n\equiv-c_1\bar{z}_1-c_2\bar{z}_2\,({\rm mod}\,c_1c_2d)}}
\frac1{c_1+c_2}.
\end{eqnarray*}
If we use $3d(c_1+c_2)+d+4.5\frac{(c_1+c_2)^2}{c_1c_2}$ as 
the number of terms in the $z_1,z_2,m,n$ sums as proved in 
Lemma \ref{congr}, we have
\begin{eqnarray}\label{T11c1c2}
T_{11}^{\eta_1\eta_2} 
& \ll & 
K_1^{1+\varepsilon}L_1L_2
\sum_{d\leq\frac{X}{K_1^{1-\varepsilon}L_1}}
\sum_{\genfrac{}{}{0pt}{1}
{c_j\leq\frac{X}{dK_j^{1-\varepsilon}L_j}}{(c_1,c_2)=1}}
\Big(1+\frac{c_1+c_2}{c_1c_2d}\Big)
\\
& \ll & 
K_1^{1+\varepsilon}L_1L_2X
\sum_{c_1\leq\frac{X}{K_1^{1-\varepsilon}L_1}}
\sum_{c_2\leq\frac{X}{K_1^{1-\varepsilon}L_2}}
1
\ll
\frac{X^{3+\varepsilon}}{K_1},
\nonumber
\end{eqnarray}
because $\frac1{c_1}+\frac1{c_2}\leq2$. 

\setcounter{equation}{0}
\section{Case of $\eta_1\neq\eta_2$ for $D_{11}$ 
when $0<\beta<1$} 

Going back to \eqref{dutheta}, we observe that for 
$\eta_1\neq\eta_2$ and $K_1=K_2$, we have 
\begin{eqnarray}
\sqrt{R_1}-\sqrt{R_2}
&=&
\frac{R_1-R_2}{\sqrt{R_1}+\sqrt{R_2}},
\label{Rquot}
\\
R_1-R_2
&=&
\frac{16\pi^2uv}{d^2}
\Big(\frac1{c_1^2}-\frac1{c_2^2}\Big)
=
\frac{16\pi^2uv}{d^2}
\Big(\frac1{c_1}-\frac1{c_2}\Big)
\Big(\frac1{c_1}+\frac1{c_2}\Big),
\label{Rdiff}
\\
\sqrt{R_1}+\sqrt{R_2}
&=&
\frac{4\pi\sqrt{uv}}d
\Big(\frac1{c_1}+\frac1{c_2}\Big)
\Big(1+
O\Big(\frac{d^2K_1^2}{X^2}(c_1^2+c_2^2)\Big)
\Big).
\label{Rsum}
\end{eqnarray}
Then \eqref{Mid2Terms} becomes
$$
\frac{|c_1-c_2|}{\sqrt2c_1c_2d}
\leq
\Big|
-\frac{\eta_1\sqrt{R_1}}{4\pi u}
-\frac{\eta_2\sqrt{R_2}}{4\pi u}
\Big|
\leq
\frac{\sqrt2|c_1-c_2|}{c_1c_2d}
$$
We will first consider the case of $c_1\neq c_2$ with 
$K_1=K_2$ and \eqref{KLlarge}. Then 
$$
\frac{|c_1-c_2|}{\sqrt2c_1c_2d}
\geq
\frac1{\sqrt2c_1c_2d^2}
\geq
\frac{K_1^{1-\varepsilon}L_1K_1^{1-\varepsilon}L_2}
{\sqrt2X^2}
$$
dominates $-\alpha\beta u^{\beta-1}\asymp X^{\beta-1}$. 
Consequently, we may take $\tau=1.5|c_1-c_2|$ and 
apply Lemma \ref{congr} to get \eqref{T11J} for this 
$\tau$. 

Now we bound $J$ in \eqref{T11J} in the case at present. 
To compute \eqref{UV2-5} we note that \eqref{R1+R2} equals
\begin{equation}\label{R1-R2Sq}
\frac{(K_1-1)^2(\sqrt{R_1}-\sqrt{R_2})^2}
{16\pi^2u^2v^2\sqrt{R_1R_2}}
=
\frac{(K_1-1)^2(R_1-R_2)^2}
{16\pi^2u^2v^2\sqrt{R_1R_2}(\sqrt{R_1}+\sqrt{R_2})^2}
\end{equation}
by \eqref{Rquot}. Using \eqref{Rdiff}, \eqref{Rsum} 
and 
$$
\sqrt{R_1R_2}
=
\frac{16\pi^2uv}{c_1c_2d^2}
\Big(1+
O\Big(\frac{d^2K_1^2}{X^2}(c_1^2+c_2^2)\Big)
\Big),
$$
\eqref{R1-R2Sq} and hence \eqref{R1+R2} and \eqref{UV2-5} 
are equal to 
$$
\frac{c_1c_2K_1^2}{16\pi^2u^2v^2}
\Big(\frac1{c_1}-\frac1{c_2}\Big)^2
\Big(1+
O\Big(\frac{d^2K_1^2}{X^2}(c_1^2+c_2^2)\Big)
\Big)
=
\frac{K_1^2|c_1-c_2|^2}{16\pi^2u^2v^2c_1c_2}
\Big(1+
O\Big(\frac{d^2K_1^2}{X^2}(c_1^2+c_2^2)\Big)
\Big).
$$

By \eqref{Rj-/Rj} and \eqref{KLlarge}, we know that 
$$
\Big|
\sum_{i=2}^5
U_i
\Big|
=
\frac1{2du}
\sqrt{\frac vu}
\Big|\frac1{c_1}-\frac1{c_2}\Big|
\Big(1+
O\Big(\frac{d^2K_1^2}{X^2}(c_1^2+c_2^2)\Big)
\Big)
\gg
X^{\beta-2+\varepsilon}
\gg
U_1\asymp V_1.
$$
Likewise, $\sum_{i=2}^5V_i$ dominates $U_1\asymp V_1$. 
Then \eqref{V1U}, \eqref{r1^4}, and \eqref{r1r2^2} 
still hold and the left hand side of \eqref{-uvSq} is 
still $\asymp$ \eqref{UV2-5B} if we replace 
$(\frac1{c_1}+\frac1{c_2})$ by $|\frac1{c_1}-\frac1{c_2}|$. 
Since \eqref{r1^4} is still a power larger than 
\eqref{r1r2^2} with $|\frac1{c_1}-\frac1{c_2}|$, we 
may choose 
$$
r_1
=
r_2
=
\frac{|c_1-c_2|^\frac12K_1^\frac12}{(c_1c_2)^\frac14X}
$$
as in \eqref{r1r2B}.
By these $r_1$ and $r_2$ and \eqref{varB}, 
Lemma \ref{2ndDer} implies that 
\begin{equation}\label{-JB}
J
\ll
\frac{{\rm var}(a)}{r_1r_2}
\ll
\frac{c_1c_2dX^{1+\varepsilon}}
{K_1|c_1-c_2|}. 
\end{equation}
Applying \eqref{-JB} to \eqref{T11J} with $\tau=1.5|c_1-c_2|$, 
the contribution of the terms with $c_1\neq c_2$ to \eqref{T11J} 
is
\begin{eqnarray}\label{c1n=c2B}
T_{c_1\neq c_2}^{\eta_1\eta_2} 
& \ll & 
K_1L_1K_2L_2
\sum_{d\leq\min(
\frac{X}{K_1^{1-\varepsilon}L_1},
\frac{X}{K_2^{1-\varepsilon}L_2}
) 
}
\frac{1}{d^2} 
\sum_{\genfrac{}{}{0pt}{1}
{c_j\leq\frac{X}{dK_j^{1-\varepsilon}L_j}}
{(c_1,c_2)=1}}
\frac{d|c_1-c_2|}{c_1c_2}
\frac{c_1c_2dX^{1+\varepsilon}}
{K_1|c_1-c_2|}
\\
&\ll& 
K_1L_1L_2X^{1+\varepsilon}
\sum_{c_1\leq\frac{X}{K_1^{1-\varepsilon}L_1}}
\sum_{c_2\leq\frac{X}{K_1^{1-\varepsilon}L_2}}
1
\ll
\frac{X^{3+\varepsilon}}{K_1},
\nonumber
\end{eqnarray}
when $\eta_1\neq\eta_2$, $K_1=K_2$, and 
$c_1\neq c_2$ under \eqref{KLlarge}. 

Now let us turn to the terms with $c_1=c_2$. Since 
$(c_1,c_2)=1$, this means that $c_1=c_2=1$. Then \eqref{dutheta} 
becomes 
$\frac{\partial\theta}{\partial u}
=-\alpha\beta u^{\beta-1}-\frac md$. If $m\neq0$, then 
$$
\Big|\frac md\Big|
\geq
\frac1d
\geq
\frac1{\min(
\frac{X}{K_1^{1-\varepsilon}L_1},
\frac{X}{K_2^{1-\varepsilon}L_2}
)}
=
\frac{K_1^{1-\varepsilon}}X
\max(L_1,L_2)
$$
dominates $-\alpha\beta u^{\beta-1}\asymp X^{\beta-1}$. 
Consequently, $\frac{\partial\theta}{\partial u}\asymp\frac md$ 
when $m\neq0$. As in \S5, $J$ is then negligible. The same is 
true for $n\neq0$. When $m=n=0$, we have 
$\frac{\partial\theta}{\partial u}\asymp X^{\beta-1}$, and hence 
we may set $T=X^\beta$ and $M=X$. We also have $U=\frac dX$ and 
$N=X$ as in \S5. Then by the first derivative test in 
\cite[Theorem 1.1]{MSY} again, the $u$-integral in $J$ is 
$\ll\frac d{X^{(n+1)\beta}}$ which is negligible for $n$ 
sufficiently large. Therefore, the contribution of the 
case of $\eta_1\neq\eta_2$ and $c_1=c_2$ is negligible. 

By \eqref{c1n=c2B} we obtain 
$T_{11}^{\eta_1\eta_2} 
\ll
\frac{X^{3+\varepsilon}}{K_1}$ 
when $\eta_1\neq\eta_2$ and $K_1=K_2$ under \eqref{KLlarge}. 
Together with the bound in \eqref{T11c1c2} for the 
case of $\eta_1=\eta_2$, we finally prove
\begin{equation}\label{D11B}
D_{11}
\ll
\frac{X^{3+\varepsilon}}{K_1},
\end{equation}
when $K_1=K_2$ under \eqref{KLlarge} when $0<\beta<1$. 

Recall that $D_{01}$ and $D_{10}$ are negligible when 
$K_1=K_2\geq X^{\frac12+\varepsilon}$ by 
\eqref{T01B} and \eqref{T10B}. When 
$K_1=K_2\leq X^{\frac12+\varepsilon}$, we have 
$\frac{X^3}{K_1}\geq X^\frac52\geq L_1L_2X^\frac32$ 
and hence the bound in \eqref{D11B} donimates 
those in \eqref{T01B} and \eqref{T10B}. 
Back to \eqref{4Terms}, 
by collecting \eqref{D00B} and \eqref{D11B} we prove 
\eqref{K1=K2>} for Theorem \ref{Thm1} when $0<\beta<1$. 

\setcounter{equation}{0}
\section{Case of $\beta=1$ for $D_{11}$}

When $\beta=1$, \eqref{dutheta} becomes 
$$
\frac{\partial\theta^{\eta_1\eta_2}}{\partial u}
=
-\alpha\beta 
-\frac{\eta_1}{4\pi u}
\sqrt{R_1}
-\frac{\eta_2}{4\pi u}
\sqrt{R_2}
-\frac{m}{c_1c_2d}.
$$
By the same arguments after \eqref{dutheta}, $J$ is 
neglible for 
$|m+\alpha c_1c_2d|\geq\tau$ and for 
$|n+\alpha c_1c_2d|\geq\tau$. Consequently, \eqref{T11J} 
holds after replacing summation conditions $|m|\leq\tau$ 
and $|n|\leq\tau$ by $|m+\alpha c_1c_2d|\leq\tau$ 
and $|n+\alpha c_1c_2d|\leq\tau$. Then Lemma \ref{congr} 
remains valid. 

To use the second derivative test, we observe that 
$U_1=V_1=0$ in \eqref{uu} and \eqref{vv}. Then we don't 
need to assume \eqref{KLlarge}, and the same calculation 
leads to \eqref{c1n=c2B} for the case of $\eta_1=\eta_2$ 
and \eqref{T11c1c2} for the case of $\eta_1\neq\eta_2$ with 
$c_1\neq c_2$. 

When $\eta_1\neq\eta_2$ with $c_1=c_2=1$, 
$\frac{\partial\theta}{\partial u}=-\alpha-\frac md$.  
Recall that 
$$
\frac{\partial^s}{\partial u^s}a(u,v)
\ll
\frac U{N^s},
\ \ \ \ 
U=\frac{dc_1^\frac12c_2^\frac12}X,
\ \ 
N=\frac X{K_1^\varepsilon}.
$$
If $|\alpha+\frac md|\geq X^{\varepsilon-1}$, then each 
integration by parts in 
$\int_{\mathbb R} a(u,v)e(\theta(u,v)) du$ produces a 
factor $\ll|\alpha+\frac md|^{-1}K_1^\varepsilon X^{-1} 
\ll X^{-\varepsilon}K_1^\varepsilon$. Hence the integral 
with respect to $u$ is 
negligible when $|\alpha+\frac md|\geq X^{\varepsilon-1}$. 
Similarly the $v$-integral is negligible for 
$|\alpha-\frac nd|\geq X^{\varepsilon-1}$. 
Then the corresponding sum for \eqref{T11afterP} becomes 
\begin{eqnarray}\label{beta=1c=1}
T_{c_1=c_2=1}^{\eta_1\neq\eta_2} 
& = & 
K_1L_1K_2L_2
\sum_{d\leq\min(\frac{X}{K_1^{1-\varepsilon}L_1}, 
\frac{X}{K_2^{1-\varepsilon}L_2})}
\frac{1}{d^2}
\sideset{}{^{*}}\sum_{z_1,z_2\,{\rm mod}\,d}
\\
& \times & 
\sum_{\genfrac{}{}{0pt}{1}
{m\equiv -z_1-z_2\,({\rm mod}\,d)}
{|\alpha+\frac md|\leq X^{\varepsilon-1}}}
\ 
\sum_{\genfrac{}{}{0pt}{1}
{n\equiv -\bar{z}_1-\bar{z}_2\,({\rm mod}\,d)}
{|\alpha-\frac nd|\leq X^{\varepsilon-1}}}
J.
\nonumber
\end{eqnarray}

Under the assumption of 
$\max(K_1L_1,K_2L_2)=X^\delta$ with $0<\delta<1$, we have 
$$
\min\Big(
\frac{X}{K_1^{1-\varepsilon}L_1}, 
\frac{X}{K_2^{1-\varepsilon}L_2}
\Big)
\leq X^{1-\delta+\varepsilon},
$$
and hence 
\begin{equation}\label{Power>1}
\Big(
\min\Big(
\frac{X}{K_1^{1-\varepsilon}L_1}, 
\frac{X}{K_2^{1-\varepsilon}L_2}
\Big)
\Big)^{\frac{1-\varepsilon}{1-\delta+\varepsilon}}
\leq 
X^{1-\varepsilon}. 
\end{equation}
Take $\varepsilon<\frac\delta2$ so that the exponent in 
\eqref{Power>1} becomes 
$\frac{1-\varepsilon}{1-\delta+\varepsilon}>1$. 
Consequently, for a given $d$ as in \eqref{beta=1c=1}, 
there is at most one $m$ and $n$ satisfying 
$|\alpha+\frac md|\leq X^{\varepsilon-1}$ and 
$|\alpha-\frac nd|\leq X^{\varepsilon-1}$. For such a 
triple $d,m,n$, taking any $z_1\ {\rm mod}\ d$ with 
$(z_1,d)=1$, there is at most one $z_1\ {\rm mod}\ d$ 
satisfying the congruences in \eqref{beta=1c=1}. Thus, 
the multiple sums in \eqref{beta=1c=1} have at most 
$d\min(
\frac{X}{K_1^{1-\varepsilon}L_1}, 
\frac{X}{K_2^{1-\varepsilon}L_2}
)$ terms. 

By $|a(u,v)|\leq\frac dX$, we bound $J\ll dX$ trivially. 
Then \eqref{beta=1c=1} is bounded by
\begin{equation}\label{beta=1c=1B}
\ll
K_1L_1K_2L_2X
\min\Big(
\frac{X}{K_1^{1-\varepsilon}L_1}, 
\frac{X}{K_2^{1-\varepsilon}L_2}
\Big)
\ll
\min(K_1L_1,K_2L_2)X^{2+\varepsilon}.
\end{equation}
Collecting \eqref{beta=1c=1B}, \eqref{c1n=c2B} for the 
case of $\eta_1=\eta_2$, and \eqref{T11c1c2} for the 
case of $\eta_1\neq\eta_2$ with $c_1\neq c_2$, we prove 
\eqref{K1=K2beta=1}.

\end{document}